\documentclass{article}

\usepackage{amsmath,amssymb}
\usepackage{url}
\newtheorem{thm}{Theorem}

\begin{document}

\title{P\'olya's Random Walk Theorem}
\author{Jonathan Novak}
\date{}
\maketitle

\begin{abstract}
	This note presents a proof of P\'olya's random walk theorem
	using classical methods from special function theory and asymptotic 
	analysis.
\end{abstract}

\section{Introduction.}

This note is about a remarkable law of nature discovered by George P\'olya \cite{Polya}.
Consider a particle situated at a given point of the integer lattice $\mathbb{Z}^d$.
Suppose that, at each tick of the clock, the particle jumps to a randomly selected neighbouring
lattice point, with equal probability of jumping in any direction.  In other words,
this particle is executing the \emph{simple random walk} on $\mathbb{Z}^d$.

A random walk is said to be \emph{recurrent} if it returns to its initial position with probability one.
A random walk which is not recurrent is called \emph{transient}.  P\'olya's classic result \cite{Polya}
is the following.

	\begin{thm}
		\label{thm:main}
		The simple random walk on $\mathbb{Z}^d$ is recurrent in dimensions $d=1,2$ and
		transient in dimension $d \geq 3$.
	\end{thm}

P\'olya's theorem is a foundational result in the theory of random walks, and many proofs are available.	
This note presents a new proof of P\'olya's Theorem using techniques
developed by de Moivre and Laplace in the eighteenth and nineteenth centuries in order to
to establish the basic limit theorems of probability theory,
see \cite[Chapter 2]{Fischer}.  These classical methods have returned to the forefront of
contemporary probability, where new universality classes of limit theorems are being investigated 
via the asymptotic analysis of exact formulas \cite{BG}.  Thus, in a sense, the proof of P\'olya's 
theorem given here is both more classical and more modern than the arguments one finds in textbooks.

\section{Loop decomposition}
Let $E$ denote the event that the simple random walk on $\mathbb{Z}^d$ returns to its initial position,
and put $p=\operatorname{Prob}(E)$.  For $n \geq 1$, let
$E_n$ be the event that the random walk returns to its initial position for the first time after $n$ 
steps.  It is convenient to set $E_0=\emptyset$, corresponding to the fact that the initial position of
the random walk does not count as a return (if it did, the return probability of any random walk would
be one).  The events $E_n$ are mutually exclusive for different values of $n$, and

	\begin{equation*}
		E = \bigsqcup_{n\geq 0} E_n.
	\end{equation*}
	
\noindent
Hence

	\begin{equation*}
		p = \sum_{n \geq 0} p_n,
	\end{equation*}
	
\noindent
where $p_n=\operatorname{Prob}(E_n)$.

A \emph{loop} on $\mathbb{Z}^d$ is a walk which begins and ends at a given point.  It is 
convenient to consider walks of length zero as loops; such loops are called \emph{trivial}.  
A non-trivial loop is \emph{indecomposable} if it is not the concatenation of two non-trivial loops.
Choose a particular point of $\mathbb{Z}^d$, and let $\ell_n$ denote the number of loops of length $n$ based at this
point.  Let $r_n$ denote the number of these which are indecomposable.  Note that $\ell_0=1$ while $r_0=0$.  
Since any non-trivial loop is the concatenation of an indecomposable loop followed by a (possibly trivial) loop, 
the counts $\ell_n$ and $r_n$ are related by

	\begin{equation*}
		\ell_n = \sum_{k=0}^n r_k \ell_{n-k}
	\end{equation*}
	
\noindent
for all $n \geq 1$.  Dividing both sides of this equation by $(2d)^n$, the total number of length $n$ walks emanating 
from a given point of $\mathbb{Z}^d$, we obtain the relation

	\begin{equation*}
		q_n = \sum_{k=0}^n p_kq_{n-k}
	\end{equation*}
	
\noindent
for all $n \geq 1$, where as above $p_n$ is the probability that the random walk returns to its 
initial position for the first time after $n$ steps, while $q_n$ is the probability that the random 
walk is located at its original position after $n$ steps.  

We introduce the generating functions

	\begin{equation*}
		P(z) = \sum_{n=0}^{\infty} p_nz^n \quad \text{and} \quad Q(z) = \sum_{n=0}^{\infty} q_nz^n.
	\end{equation*}
	
\noindent
The relation between $p_n$ and $q_n$ is then equivalent to the identity

	\begin{equation*}
		P(z)Q(z) = Q(z)-1
	\end{equation*}
	
\noindent
in the algebra $\mathbb{Q}[[z]]$ of formal power series.  
Since $p_n \leq q_n \leq 1$, each of these series has radius of convergence at least one,
and the above may be considered as an identity in the algebra of analytic functions 
on the open unit disc in $\mathbb{C}$.  The function $Q(z)$ is non-vanishing for $z$ in the interval $[0,1)$,
and hence we have

	\begin{equation*}
		P(z) = 1-\frac{1}{Q(z)}, \quad z \in [0,1).
	\end{equation*}
	
\noindent
Since 

	\begin{equation*}
		P(1) = \sum_{n=0}^{\infty} p_n = p,
	\end{equation*}
	
\noindent
Abel's power series theorem applies and we have

	\begin{equation*}
		p = \lim_{\substack{z \rightarrow 1\\ z\in [0,1)}} P(z) = 1 - \frac{1}{\lim\limits_{\substack{z \rightarrow 1\\ z\in [0,1)}} Q(z)}.
	\end{equation*}
	
\noindent
The limit in the denominator is either $+\infty$ or a positive real number.  In the former case we have $p=1$ (recurrence),
and in the latter $p<1$ (transience).

\section{Exponential loop generating function}
In order to analyze the limit in question, we need
a tractable representation of the function $Q(z)$.  This amounts to finding an 
expression for the loop generating function 

	\begin{equation*}
		L(z) = \sum_{n=0}^{\infty} \ell_n z^n.
	\end{equation*}
	
\noindent
Indeed, $Q(z) = L(\frac{z}{2d})$.

While the ordinary generating function $L(z)$ is difficult to analyze directly, the exponential
loop generating function

	\begin{equation*}
		E(z) = \sum_{n=0}^{\infty} \ell_n \frac{z^n}{n!}
	\end{equation*}
	
\noindent
is quite accessible.  This is because any loop on $\mathbb{Z}^d$ is a shuffle of loops on 
$\mathbb{Z}^1$, and products of exponential generating functions correspond to shuffles.
This is a basic property of exponential generating functions which we will review in the
specific case at hand.  For a general treatment, the reader is referred to \cite[Chapter 5]{Stanley}.
 
In this paragraph it is important to make the dependence on $d$ explicit, 
so we write $\ell_n^{(d)}$ for the number of length $n$ loops on $\mathbb{Z}^d$ and 
$E_d(z)$ for the exponential generating function of this sequence.
Let us consider the case $d=2$.  A loop on $\mathbb{Z}^2$ is a closed walk
which takes unit steps in two directions, horizontal and vertical.  A length $n$ loop
on $\mathbb{Z}^2$ is made up of some number $k$ of horizontal steps together with $n-k$ vertical steps.  
The $k$ horizontal steps constitute a length $k$ loop on $\mathbb{Z}$, and the $n-k$ vertical steps 
constitute a length $n-k$ loop on $\mathbb{Z}$.  Thus, the number of length $n$ loops on 
$\mathbb{Z}^2$ which take $k$ horizontal and $n-k$ vertical steps is 

	\begin{equation*}
		{n \choose k} \ell_k^{(1)} \ell_{n-k}^{(1)},
	\end{equation*}
	
\noindent
since specifying the times at which the $k$ horizontal steps occur 
uniquely determines the times at which the $n-k$ vertical steps occur.  
The total number of length $n$ loops on $\mathbb{Z}^2$ is therefore

	\begin{equation*}
		\ell_n^{(2)} = \sum_{k=0}^{n} {n \choose k} \ell_k^{(1)} \ell_{n-k}^{(1)}.
	\end{equation*}
	
\noindent
This is equivalent to the generating function identity

	\begin{equation*}
		E_2(z) = E_1(z)^2.
	\end{equation*}
	
\noindent
The same reasoning applies for any $d$, and in general we have

	\begin{equation*}
		E_d(z) = E_1(z)^d.
	\end{equation*}
	
Counting loops in one dimension is easy,

	\begin{equation*}
		\ell_n^{(1)} = \begin{cases}
			{2k \choose k}, \text{ if } n=2k \text{ is even,} \\
			0, \text{ if } n \text{ is odd.}
			\end{cases}.
	\end{equation*}
	
\noindent
Indeed, any loop on $\mathbb{Z}$ consists of $k$ positive steps and $k$ negative steps for some $k \geq 0$,
and the times at which the positive steps occur determine the times at which the negative steps occur.
Thus

	\begin{equation*}
		E_1(z) = \sum_{k=0}^\infty {2k \choose k} \frac{z^{2k}}{(2k)!} = \sum_{k=0}^\infty \frac{z^{2k}}{k!k!}.
	\end{equation*}
	
\noindent
Now a minor miracle occurs: the exponential generating function for lattice walks in one dimension
is a \emph{modified Bessel function of the first kind}.

The modified Bessel function of the first kind, usually denoted $I_\alpha(z)$, 
is one of two linearly independent solutions to the second order differential equation

	\begin{equation*}
		\bigg{(} z^2 \frac{d^2}{dz^2} + z \frac{d}{dz} - (z^2+\alpha^2)\bigg{)}F(z) =0, \quad \alpha \in \mathbb{C}.
	\end{equation*}
	
\noindent
This differential equation is known as the \emph{modified Bessel equation}; it appears in a multitude of physical
problems, and was exhaustively studied by nineteenth century mathematicians.  An excellent reference on this subject
is \cite[Chapter 4]{AAR}.  It is known that the modified Bessel function admits both a series representation,

	\begin{equation*}
		I_\alpha(z) = \sum_{k=0}^{\infty} \frac{(\frac{z}{2})^{2k+\alpha}}{k!\Gamma(k+\alpha+1)},
	\end{equation*}
	
\noindent
and an integral representation,

	\begin{equation*}
		I_\alpha(z) = \frac{(\frac{z}{2})^\alpha}{\sqrt{\pi}\Gamma(\alpha+\frac{1}{2})}\int_0^\pi 
		e^{(\cos \theta)z}(\sin \theta)^{2\alpha} d\theta.
	\end{equation*}
	
\noindent
From the series representation, we see that $E_1(z) = I_0(2z)$, and hence

	\begin{equation*}
		E(z) = I_0(2z)^d.
	\end{equation*}

\section{Borel transform}
We now have a representation of the exponential generating function $E(z)$ counting loops on $\mathbb{Z}^d$
in terms of a standard mathematical object, the modified Bessel function $I_0(z)$.  What we need, however,
is a representation of the ordinary loop generating function $L(z)$.  

The integral transform

	\begin{equation*}
		(\mathcal{B}f)(z) = \int_0^{\infty} f(tz)e^{-t} dt,
	\end{equation*}
	
\noindent
which looks like the Laplace transform of $f$ but with the $z$-parameter in the wrong place,
converts exponential generating functions into ordinary generating functions.  To see why,
write out the Maclaurin series of $f(tz)$, interchange integration and summation to obtain

	\begin{equation*}
		(\mathcal{B}f)(z) = \sum_{n=0}^{\infty} f^{(n)}(0) \frac{z^n}{n!} \int_0^{\infty} t^ne^{-t} dt,
	\end{equation*}
	
\noindent
and use the fact that

	\begin{equation*}
		\int_0^{\infty} t^ne^{-t} dt = n!.
	\end{equation*}

\noindent
The transform $f\mapsto \mathcal{B}f$ was invented 
by Borel in order to ``sum'' divergent series \cite[p. 55]{Borel}.
In our case, the Borel transform produces the formula

	\begin{equation*}
		L(z) = \mathcal{B}E(z) = \mathcal{B}I_0(2z)^d = \int_0^{\infty} I_0(2tz)^d e^{-t} dt,
	\end{equation*}
	
\noindent
which in turn leads to the integral representation 

	\begin{equation*}
		Q(z) = L(\frac{z}{2d}) = \int_0^{\infty} I_0\bigg{(}\frac{tz}{d}\bigg{)}^d e^{-t} dt.
	\end{equation*}

\section{The Laplace principle}
We will now use the integral representation just obtained to determine whether the 
limit under consideration is finite or infinite.  It suffices to answer this question
for the tail integral

	\begin{equation*}
		\int_N^{\infty} I_0\bigg{(}\frac{tz}{d}\bigg{)}^d e^{-t} dt, \quad N \gg 0.
	\end{equation*}
	
\noindent
For $N$ large, the behaviour of the tail integral is in turn determined by the 
behaviour of the integrand as $t \rightarrow \infty$.
In order to estimate the integrand, we invoke the formula

	\begin{equation*}
		I_0\bigg{(} \frac{tz}{d} \bigg{)} = \frac{1}{\pi} \int_0^\pi e^{tf(\theta)} d\theta,
	\end{equation*}
	
\noindent
where $f(\theta) = \frac{z}{d} \cos \theta$,
and estimate this integral as $t \rightarrow \infty$ using a basic technique of asymptotic
analysis known as \emph{Laplace's method}.  

The function $f(\theta)$ is strictly maximized
over the interval $[0,\pi]$ at the left endpoint $\theta=0$.  Thus the integrand 
$e^{tf(\theta)}$ is exponentially larger at $\theta=0$ than at any other point of this interval.  As
$t \rightarrow \infty$ this effect becomes increasingly exaggerated, so much so that the 
integral ``localizes'' at $\theta=a$ in the $t \rightarrow \infty$ limit.  To quantify this, note 
that $f'(0)=0, f''(0)<0$, and consider the quadratic Taylor approximation of $f(\theta)$:

	\begin{equation*}
		f(\theta) \approx f(0) - |f''(0)|\frac{\theta^2}{2}.
	\end{equation*}
	
\noindent 
Replacing $f(\theta)$ with its quadratic approximation, we obtain the integral
approximation

	\begin{equation*}
		\int_0^\pi e^{tf(\theta)} d\theta \approx e^{tf(0)} \int_0^\pi e^{-t|f''(0)|\frac{\theta^2}{2}} d\theta.
	\end{equation*}	
	
\noindent
Extending the integral on the right over the positive reals and ignoring the rapidly decaying error 
incurred results in a half a Gaussian integral, which can be computed exactly:

	\begin{equation*}
		\int_0^{+\infty} e^{-t|f''(0)|\frac{\theta^2}{2}} d\theta = \sqrt{\frac{\pi}{2t |f''(0)|}}.
	\end{equation*}
	
\noindent
Thus we expect that 

	\begin{equation*}
		\int_0^\pi e^{tf(\theta)} d\theta \approx e^{tf(0)} \sqrt{\frac{\pi}{2t |f''(0)|}}
	\end{equation*}
	
\noindent
is an approximation of our integral whose accuracy increases as $t \rightarrow \infty$.  Laplace's 
principle (see e.g. \cite[\S 5.2]{Bleistein}) is the statement that this is indeed the case;  we have

	\begin{equation*}
		\int_0^\pi e^{tf(\theta)} d\theta \sim e^{tf(0)} \sqrt{\frac{\pi}{2t |f''(0)|}}, \quad t \rightarrow \infty,
	\end{equation*}
	
\noindent
where the notation $F(t) \sim G(t),\ t \rightarrow \infty$ means that $\lim_{t \rightarrow \infty} \frac{F(t)}{G(t)}=1$.
	
Putting everything together, we have the asymptotic formula

	\begin{equation*}
		I_0\bigg{(} \frac{tz}{d} \bigg{)}^d e^{-t} \sim \text{constant} \cdot e^{t(z-1)} (tz)^{-\frac{d}{2}}, \quad t \rightarrow \infty.
	\end{equation*}

\noindent
Applying the monotone convergence theorem, we find

	\begin{equation*}
		\lim_{\substack{z\rightarrow 1\\ z\in [0,1)}} \int_N^{\infty} e^{t(z-1)}(tz)^{-\frac{d}{2}} dt
		=  \int_N^{\infty} \lim_{\substack{z\rightarrow 1\\ z\in [0,1)}} e^{t(z-1)}(tz)^{-\frac{d}{2}} dt
		= \int_N^{\infty} t^{-\frac{d}{2}} dt,
	\end{equation*}
	
\noindent
and thus conclude that the recurrence or transience 
of the simple random walk on $\mathbb{Z}^d$ is equivalent to the divergence or 
convergence of the integral

	\begin{equation*}
		\int_N^{\infty} t^{-d/2} dt, \quad N \gg 0.
	\end{equation*}
	
\noindent
Since this integral diverges for $d=1,2$ and converges for $d \geq 3$,
P\'olya's Theorem is proved.

\bigskip

\noindent\textit{Department of Mathematics,
Massachusetts Institute of Technology, Cambridge, MA 02139\\
jnovak@math.mit.edu}

\end{document}